\documentclass[12pt]{amsart}

\usepackage{color}
\usepackage{amsmath, amsthm, amssymb}
\usepackage{amsfonts}
\usepackage[ansinew]{inputenc}
\usepackage[dvips]{epsfig}
\usepackage{graphicx}
\usepackage[english]{babel}
\usepackage{hyperref}
\theoremstyle{plain}
\newtheorem{thm}{Theorem}[section]

\theoremstyle{definition}

\theoremstyle{remark}

\numberwithin{equation}{section}

\begin{document}

\title[Euclidean balls solve isoperimetric problems with nonradial weights]{Euclidean balls solve some isoperimetric problems with nonradial weights}

\author{Xavier Cabr\'e}

\address{ICREA and Universitat Polit\`ecnica de Catalunya, Departament de Matem\`{a}tica  Aplicada I, Diagonal 647, 08028 Barcelona, Spain}
\email{xavier.cabre@upc.edu}

\author{Xavier Ros-Oton}

\address{Universitat Polit\`ecnica de Catalunya, Departament de Matem\`{a}tica  Aplicada I, Diagonal 647, 08028 Barcelona, Spain}
\email{xavier.ros.oton@upc.edu}

\thanks{
The authors were supported by grants MTM2008-06349-C03-01, MTM2011-27739-C04-01 (Spain) and 2009SGR345 (Catalunya).}

\author{Joaquim Serra}

\address{Universitat Polit\`ecnica de Catalunya, Departament de Matem\`{a}tica  Aplicada I, Diagonal 647, 08028 Barcelona, Spain}

\email{joaquim.serra@upc.edu}

\maketitle

\begin{abstract}

\vskip 0.5\baselineskip

In this note we present the solution of some isoperimetric problems in open convex cones of $\mathbb R^n$ in which perimeter and volume are measured with respect to certain nonradial weights.
Surprisingly, Euclidean balls centered at the origin (intersected with the convex cone) minimize the isoperimetric quotient.
Our result applies to all nonnegative homogeneous weights satisfying a concavity condition in the cone.
When the weight is constant, the result was established by Lions and Pacella in 1990.

\vskip 1\baselineskip

\noindent{\sc R\'esum\'e.}
Dans cette note, nous pr\'esentons la solution de certains probl\`emes isop\'erim\'etriques dans des c\^ones convexes de $\mathbb R^n$  o\`u le p\'erim\`etre et le volume sont mesur\'es par rapport \`a certains poids non radiaux.
Contrairement \`a ce que l'on pourrait penser, les boules euclidiennes centr\'ees \`a l'origine (intersect\'ees avec le c\^one) minimisent le quotient isop\'erim\'etrique.
Notre r\'esultat s'applique aux poids strictement positifs, homog\`enes et satisfaisant une condition de concavit\'e dans le c\^one.
Lorsque le poids est constant, le r\'esultat a \'et\'e \'etabli par Lions et Pacella en 1990.
\end{abstract}

\vspace{6mm}

\section{Introduction}
\label{intro}

The aim of this note is to present some new isoperimetric inequalities with weights (also called densities) in convex cones of $\mathbb R^n$.
More general results will appear in \cite{CRS}.
The isoperimetric problem with a density $w$ concerns the existence and characterization of minimizers of the weighted perimeter $\int_{\partial\Omega}w$ among those sets $\Omega$ having constant weighted volume $\int_\Omega w$.
These type of problems have attracted attention recently ---see for example the survey \cite{M} in the Notices of the AMS.
However, weighted isoperimetric inequalities with best constant are known in very few cases, even in the case of the plane ($n=2$).

Here we present the solution of the isoperimetric problem in any open convex cone of $\mathbb R^n$ for certain nonradial weights.
Namely, our result applies to all nonnegative homogeneous weights satisfying a concavity condition in the cone.
A surprising fact is that Euclidean balls centered at the origin (intersected with the cone) minimize the isoperimetric quotient with these nonradial weights.

Our result, stated below, extends the isoperimetric inequality in convex cones of P. L. Lions and F. Pacella \cite{LP}.
Their result states that among all sets $\Omega$ with fixed volume contained in an open convex cone $\Sigma$, the unit ball intersected with the cone minimizes the perimeter inside the cone (that is, not counting the part of the boundary of $\Omega$ that lies on the boundary of the cone).
Our extension allows any homogeneous weight $w$ satisfying a concavity condition.
Our approach is completely different from the one in \cite{LP} and, therefore, by setting $w\equiv1$ we provide with a new proof of the theorem of Lions and Pacella.

Let $\Sigma$ be an open convex cone in $\mathbb R^n$ and denote by
\[P_{\Sigma}(\Omega):=\int_{\Sigma\cap\partial\Omega}w(x)d\sigma\qquad\textrm{and}\qquad m(\Omega)=\int_{\Omega}w(x)dx\]
the weighted perimeter and volume of $\Omega\subset\Sigma$ inside the cone $\Sigma$.

\begin{thm}\label{cones}
Let $\Sigma$ be an open convex cone in $\mathbb R^n$, and let $w$ be a continuous function in $\overline\Sigma$, positive and $C^{1,\gamma}$ in $\Sigma$ for some $\gamma\in(0,1)$, homogeneous of degree $\alpha\geq0$, and such that $w^{1/\alpha}$ is concave in case $\alpha>0$.
Then, for every Lipschitz domain $\Omega\subset\Sigma$,
\begin{equation}\label{dos}
\frac{P_{\Sigma}(\Omega) }{m(\Omega)^{\frac{D-1}{D}} }\geq  \frac{P_{\Sigma}(\Sigma\cap B_1)}{m(\Sigma\cap B_1)^{\frac{D-1}{D}}},
\end{equation}
where $D=n+\alpha$ and $B_1=B_1(0)$ is the Euclidean unit ball of $\mathbb R^n$.
\end{thm}

Note the following surprising fact:
even that the weights considered here are not radial (unless $w\equiv1$), still Euclidean balls centered at the origin (intersected with the cone) minimize this isoperimetric quotient.
Note also that we allow $w$ to vanish somewhere on $\partial\Sigma$.

The exponent $D=n+\alpha$ has a dimension flavor and can be found by a scaling argument thanks to the homogeneity of the weight.
The interpretation of $D$ as a dimension is more clear in the following example.
The monomial weights
\begin{equation}\label{monomial}
\begin{split}
&\hspace{12mm} w(x)=x_1^{A_1}\cdots x_n^{A_n} \qquad \textrm{in the cone}\\
& \Sigma=\{x\in\mathbb R^n\,:\, x_i>0\textrm{ for all }i\textrm{ such that }A_i>0\},
\end{split}
\end{equation}
where $A_i\geq0$, $\alpha=A_1+\cdots+A_n$, and $D=n+A_1+\cdots+A_n$, are important examples for which (\ref{dos}) holds.
These isoperimetric (and the corresponding Sobolev) inequalities with monomial weights are studied by the first two authors in \cite{CR}.
They arose in \cite{CR2}, where we studied reaction-diffusion problems with symmetry of double revolution.
A function $u$ has symmetry of double revolution when $u(x,y)= u(|x|,|y|)$, with  $(x,y) \in \mathbb R^{D}= \mathbb R^{A_1+1}\times\mathbb R^{A_2+1}$ (here we assume $A_i$ to be positive integers). In this way, $u$ can be seen as a function in $\mathbb R^2=\mathbb R^n$, and it is here where the Jacobian $x^{A_1}y^{A_2}$ appears.
A similar argument under multiple revolutions shows that, when $w$ and $\Sigma$ are given by (\ref{monomial}) and all $A_i$ are nonnegative integers, Theorem \ref{cones} follows from the classical isoperimetric inequality in $\mathbb R^D$.

We know only of two results where nonradial weights lead to radial minimizers.
The first one is the isoperimetric inequality by Maderna-Salsa \cite{MS} in the upper half plane $\{(x,y)\in\mathbb R^2\,:\, y>0\}$ with the weight $y^k$, $k>0$.
The second one is due to Brock-Chiacchio-Mercaldo \cite{BCM} and extends the one in \cite{MS} by including the weights $(x_n)^{k}\exp(c|x|^2)$ in $\mathbb R^n_+$, with $k\geq0$ and $c\geq0$. They prove that half balls are the minimizers of the isoperimetric quotient with these weights.
Our result provides a much wider class of nonradial weights and cones for which the optimizers are Euclidean balls.

Of course, not all homogeneous weights lead to radial minimizers. In fact, the following example shows that even radial homogeneous weights may lead to nonradial minimizers.
Indeed, consider in $\mathbb R^2$ the weight $w(x)=|x|^\alpha$, where $\alpha>0$, $|\cdot|$ is the Euclidean norm, and $\Sigma$ is an open convex cone of angle $\beta$.
Then, it is proved in \cite{DHHT} that there exists $\beta_0\in (0,\pi)$ such that, for $\beta<\beta_0$ the set $\Sigma\cap B_1$ minimizes the isoperimetric quotient in $\Sigma$ with weight $w$, while for $\beta>\beta_0$ it does not.

\section{Sketch of the proof}

The proof of Theorem \ref{cones} follows the ideas introduced by the first author in a new proof of the classical isoperimetric inequality; see \cite{CSCM,CDCDS}.
It is quite surprising (and fortunate) that this proof (which gives the best constant) can be adapted to the case of the previous homogeneous weights.

We next sketch the proof of our result.
To simplify the exposition, here we assume that $w\equiv0$ on $\partial\Sigma$ (at the end of the Note we explain what must be changed in the general case).
For example, the monomial weights (\ref{monomial}) satisfy the assumption $w\equiv0$ on $\partial\Sigma$.
As a consequence, $P_\Sigma(\Omega)=\int_{\Sigma\cap\partial\Omega}w=\int_{\partial\Omega}w:=P(\Omega)$.
Thus, by regularizing $\Omega$, we may assume that $\overline\Omega\subset\Sigma$ and $\Omega$ is smooth.

Consider the solution $u$ of the linear Neumann problem
\begin{equation} \label{eqsem}
\left\{ \begin{array}{ll} w^{-1}\textrm{div}(w\nabla u) = b_\Omega  &\quad \mbox{in } \Omega\\
\frac{\partial u}{\partial\nu} =1 &\quad \textrm{on }\partial \Omega.
\end{array}\right. \end{equation}
Note that, since $w>0$ in $\Sigma$, (\ref{eqsem}) is a uniformly elliptic equation.
Since $w\in C^{1,\gamma}(\Sigma)$, $u\in C^{2,\gamma}(\overline\Omega)$.
The constant $b$ is chosen so that (\ref{eqsem}) has a unique solution up to an additive constant, that is,
\begin{equation}\label{cttb} b_\Omega=\frac{P(\Omega)}{m(\Omega)}.\end{equation}

We now consider the lower contact set of $u$, $\Gamma_u$, defined as the set of points in $\Omega$ at which the tangent hyperplane to the graph of $u$ lies below $u$ in all $\overline \Omega$.
Then, as in the ABP method, we touch by below the graph of $u$ with hyperplanes of fixed slope $p\in B_1$, and using the boundary condition in (\ref{eqsem}) we deduce that
\[B_1 \subset \nabla u (\Gamma_u).\]
From this, we obtain $\Sigma\cap B_1\subset \Sigma\cap \nabla u(\Gamma_u)$ and thus
\begin{equation}\label{ineq}
\begin{split}
m(\Sigma\cap B_1) &\leq \int_{\Sigma\cap \nabla u (\Gamma_u)}w(p)dp \\
&\leq \int_{(\nabla u)^{-1}(\Sigma)\cap \Gamma_u} w(\nabla u(x))\det D^2u(x)\,dx\\
&\leq  \int_{(\nabla u)^{-1}(\Sigma)\cap \Gamma_u} w(\nabla u)\left(\frac{\Delta u}{n}\right)^ndx.
\end{split}
\end{equation}
We have applied the area formula to the map $\nabla u : \Gamma_u \rightarrow \mathbb R^n$ and also the classical arithmetic-geometric mean inequality ---all eigenvalues of $D^2u$ are nonnegative in $\Gamma_u$ by definition of this set.

Next we use that, when $\alpha>0$,
\[Y^{\alpha}Z^n\leq \left(\frac{\alpha Y+nZ}{\alpha+n}\right)^{\alpha+n}\ \ \textrm{for all positive}\ \ Y\ \textrm{and}\ Z,\]
and also that
\[\alpha\left(\frac{w(p)}{w(x)}\right)^{1/\alpha}\leq \frac{\nabla w(x)\cdot p}{w(x)}
\ \ \textrm{for all}\ x\ \textrm{and}\ p\ \textrm{in}\ \Sigma,\]
which is equivalent to the concavity of $w^{1/\alpha}$ (given that $w^{1/\alpha}$ is homogeneous of degree 1).
We find
\begin{equation}\label{2}
\frac{w(\nabla u)}{w(x)}\left(\frac{\Delta u}{n}\right)^n\leq
\left(\frac{\alpha\left(\frac{w(\nabla u)}{w(x)}\right)^{1/\alpha}+\Delta u}{\alpha+n}\right)^{\alpha+n}\leq
\left(\frac{\frac{\nabla w(x)\cdot \nabla u}{w(x)}+\Delta u}{\alpha+n}\right)^{\alpha+n}=
\left(\frac{b_\Omega}{D}\right)^{D}.\end{equation}
In the last equality we have used equation (\ref{eqsem}). If $\alpha=0$ then $w\equiv1$, and (\ref{2}) is trivial.

Therefore, since $\Gamma_u\subset \Omega$, combining (\ref{ineq}) and (\ref{2}) we obtain
\begin{equation}\label{7}\begin{split}
m(\Sigma\cap B_1) &\leq \int_{\Gamma_u} \left(\frac{b_\Omega}{D}\right)^{D}w(x)dx=
\left(\frac{b_\Omega}{D}\right)^{D}m(\Gamma_u)\\
&\leq \left(\frac{b_\Omega}{D}\right)^{D}m(\Omega)
= D^{-D}\frac{P(\Omega)^{D}}{m(\Omega)^{D-1}}.\end{split}\end{equation}
In the last equality we have used the value of the constant $b_\Omega$, given by (\ref{cttb}).

Finally, when $\Omega=\Sigma\cap B_1$ we consider $u(x)=|x|^2/2$ and $\Gamma_u=\Sigma\cap B_1$. Now, $u_\nu=1$ is only satisfied on $\Sigma\cap \partial\Omega$ but, since $w\equiv 0$ on $\partial \Sigma\cap \partial\Omega$, we have $b_{\Sigma\cap B_1} = P(\Sigma\cap B_1)/m(\Sigma\cap B_1)$ --- as in  (\ref{cttb}). For these concrete $\Omega$ and $u$ one verifies that all inequalities in (\ref{ineq}),(\ref{2}),(\ref{7}) are equalities. Thus, (\ref{dos}) follows.
\qed

Without the assumption $w\equiv0$ on $\partial\Sigma$, the proof is more involved. We need to consider a Neumann condition of the form $\partial u/\partial\nu=H(\nu)$ on $\partial\Omega$, with $H(\nu)=0$ in all normal directions $\nu$ to the cone and $H(\nu)=1$ for all $\nu$ which are directions interior to the cone.

\end{document}